\documentclass[a4paper,12pt]{article}
\usepackage{amssymb}

\newtheorem{teorema}{Theorem}
\newtheorem{propozitia}{Proposition}

\newtheorem{afirmatia}{Claim}
\newtheorem{corolarul}{Corollary}
\newtheorem{observatia}{Remark}

\newcommand{\R}{\mathbb R}

\begin{document}

\author{Ravi P. Agarwal\\
\textit{Department of Mathematical Sciences,}\\ \textit{Florida Institute of Technology,}\\
\textit{Melbourne, FL 32901, USA}\\
\textit{e-mail: agarwal@fit.edu}\\
Octavian G. Mustafa\\
\textit{Faculty of Mathematics, D.A.L.}\\
\textit{University of Craiova, Romania}\\
\textit{e-mail: octaviangenghiz@yahoo.com}}
\title{On a local theory of asymptotic integration for nonlinear differential equations}
\date{}
\maketitle

\begin{center}
\textit{Dedicated to the memory of Professor Cezar Avramescu}
\end{center}

\vspace{1cm}

\noindent\textbf{Abstract} By revisiting an asymptotic integration theory of nonlinear ordinary differential equations due to J.K. Hale and N.
Onuchic [Contributions Differential Equations 2 (1963), 61--75], we improve and generalize several recent results in the literature. As an
application, we study the existence of bounded positive solutions to a large class of semi-linear elliptic partial differential equations via
the subsolution-supersolution approach.

\section{Introduction}

Consider the ordinary differential equation in the superlinear case ($\lambda>1$)
\begin{eqnarray}
x^{\prime\prime}+q(t)x^{\lambda}=0,\qquad t\geq t_0\geq1,\label{ecEF1}
\end{eqnarray}
where the functional coefficient $q:[t_0,+\infty)\rightarrow[0,+\infty)$ is assumed continuous and with (eventual) isolated zeros, and
$x^{\lambda}=\vert x\vert^{\varepsilon}x$ for $\lambda=1+\varepsilon$. The equation (\ref{ecEF1}) and its variants with modified functional
arguments are frequently called the \textit{Emden-Fowler}, the \textit{Lane-Emden} or the \textit{Thomas-Fermi} equations, as a tribute to a
series of their fundamental particular cases, cf. \cite{Graef_etal1984}. Details on the analysis of such differential equations can be found in
the monographs \cite{AgarwalGraceORegan,Bellman_carte,KiguradzeChanturia_2}.

In 1955, F.V. Atkinson published a spectacular result on the oscillation of equation (\ref{ecEF1}) with a proof relying on asymptotic
integration theory via the Picard iterations.

\begin{teorema}\label{teor1Atkinson}\emph{(\cite[Theorem 1]{Atkinson_3})}
A necessary and sufficient condition for the oscillation of (\ref{ecEF1}) is given by
\begin{eqnarray}
\int^{+\infty}tq(t)dt=+\infty.\label{Atkinson_oscilcond}
\end{eqnarray}
\end{teorema}

The sufficiency part of Theorem \ref{teor1Atkinson} consisted in a demonstration of the next

\begin{afirmatia}
If Atkinson's hypothesis (\ref{Atkinson_oscilcond}) does not hold -- meaning that we can take (\cite[p. 646]{Atkinson_3})
\begin{eqnarray}
\eta=\lambda\int_{t_0}^{+\infty}tq(t)dt<1,\label{Atkinson_hyp1}
\end{eqnarray}
then the boundary valued problem
\begin{eqnarray}
\left\{\begin{array}{ll} x^{\prime\prime}+q(t)x^{\lambda}=0,\qquad
t\geq t_0,\\
\lim\limits_{t\rightarrow+\infty}x(t)=1,\qquad\lim\limits_{t\rightarrow+\infty}x^{\prime}(t)=0
\end{array}\right.\label{bvp_Atkinson}
\end{eqnarray}
has at least one solution.
\end{afirmatia}

Recently, Theorem \ref{teor1Atkinson} has been considered by several investigators \cite{DubeMingarelli_HaleOnuchic,Ehrnstrom,MustafaNA2005}. We
would like to improve hereafter upon the conclusions of \cite{MustafaNA2005} and of some other recent work. Important continuations of
Atkinson's theorem are given in \cite{MooreNehari_3,Waltman,Wong}.

To emphasize the significance of Atkinson's contribution, let us just say that, via Theorem \ref{teor1Atkinson}, \textit{the superlinear case of
the second order ordinary differential equations with positive functional coefficient $q(t)$ has a necessary and sufficient condition of
oscillation, a feature without correspondence to the linear case}, cf. \cite[p. 643]{Atkinson_3}, \cite[p. 128]{CoffWong}.

We shall discuss a different proof of Theorem \ref{teor1Atkinson} by controlling the behavior of the derivative $x^{\prime}(t)$ for the
solutions of (\ref{ecEF1}). This type of analysis has been already performed in different circumstances by Coffman and Wong \cite[Theorems 1,
3]{CoffWong}, and also in a general setting by Hale and Onuchic \cite{HaleOnuchic}.

Our main motivation in this work comes from a problem concerning the existence of positive classical solutions to the semi-linear elliptic partial differential equation
\begin{eqnarray}
\Delta u+f(x,u)+g(\vert x\vert)x\cdot \nabla u=0,\qquad x\in G_A,  \label{main_eq}
\end{eqnarray}
where $G_A=\{x\in {\R}^n:\vert x\vert>A\}$ and $n\geq 3$. One can find different conclusions regarding (\ref{main_eq}) in the contributions
\cite{AgarwalMustafa,Agarwal_et_al,Constantin1997,Constantin2005,Deng,MustafaPDEsJMAA}, \cite{Ehrnstrom2}--\cite{EhrnstromMustafa},
\cite{Hesaaraki}, \cite{MustafaRogovPDEs}--\cite{Orpel}, \cite{Yin}
 and in their references. In this literature, besides the regularity requirements imposed upon $f,g$, see \cite{N,NS1980}, and several technical
 conditions on $f$, the main hypothesis concerned with $g$ reads as (\cite{Constantin1997,Ehrnstrom})
\begin{eqnarray}
\int_A^{+\infty }r\vert g(r)\vert dr<+\infty .  \label{integr_1}
\end{eqnarray}

It has been noticed by Ehrnstr\"{o}m \cite{Ehrnstrom} that, when $g$ takes only nonnegative values, we can eliminate (\ref{integr_1}) provided
that a certain elliptic partial differential equation possesses a positive radially symmetric solution $U$ verifying the condition
\begin{eqnarray}
x\cdot\nabla U(x)\leq0,\qquad x\in G_B,\label{Ehrnstrom_cond}
\end{eqnarray}
where $B>A$ is taken large enough. This solution $U$ plays the role of supersolution in the comparison method approach \cite{GilbargTrudinger}
designed for the class of equations (\ref{main_eq}), cf. \cite{Constantin1997,NS1980}. Further developments of Ehrnstr\"{o}m's idea are given in
\cite{AgarwalMustafa,Agarwal_et_al,MustafaPDEsJMAA,EhrnstromMustafa}.

Translated into the language of ordinary differential equations, the research about $U$ reads as follows: find (if any) a positive solution
$x(t)$ of the equation (\ref{ecEF1}) such that $x(t)=O(t)$ as $t\rightarrow+\infty$ and
\begin{eqnarray}
x^{\prime}(t)-\frac{x(t)}{t}<0,\qquad t\geq t_0.\label{wronsk_sign}
\end{eqnarray}

Formula (\ref{wronsk_sign}) prompts the need to estimate the behavior of $x^{\prime}(t)$ for the solution $x(t)$ to (\ref{ecEF1}) when $t$
increases indefinitely and thus motivates our investigation.

This paper is organized as follows. In the next Section a presentation of the results due to Hale and Onuchic \cite{HaleOnuchic} is given and
commented upon. The hypotheses and proofs contain several modifications with respect to the original analysis. Further, these results are
specialized to the class of equations (\ref{ecEF1}) and certain estimates recently established are improved. In this process, several
theorems scattered throughout the literature are presented in a unifying way. The last Section is devoted to showing that the equation
(\ref{main_eq}) has a bounded positive solution with prescribed upper and lower bounds.

\section{The Hale-Onuchic theory of asymptotic integration}

Let us consider the boundary value problem given by the nonlinear differential equation ($k+r\geq1$)
\begin{eqnarray}
u^{(k+r+1)}=P(t,u,u^{\prime},\cdots,u^{(k+r)}),\qquad t\geq t_0\geq1,\label{ecHO}
\end{eqnarray}
where the function $P:[t_0,+\infty)\times{\R}^{k+r+1}\rightarrow{\R}$ is continuous, together with the data
\begin{eqnarray}
\left\{
\begin{array}{ll}
u^{(i)}(t_0)=a_{k-i},\qquad i\in\overline{0,k-1},\quad(\mbox{when }k\geq1)\\
\lim\limits_{t\rightarrow+\infty}u^{(k)}(t)=c,\\
\lim\limits_{t\rightarrow+\infty}[t^{j}u^{(k+j)}(t)]=0,\qquad j\in\overline{1,r}\quad(\mbox{when }r\geq1)
\end{array}
\right.\label{dateHO}
\end{eqnarray}
for fixed $a_{k-i},c\in{\R}$.

According to \cite[p. 62]{HaleOnuchic}, we can establish an existence result for the problem (\ref{ecHO}), (\ref{dateHO}) \textit{by
transforming it into an integro-differential $(r+1)$--dimensional system}.

If $u(t)$ is a solution of the problem (\ref{ecHO}), (\ref{dateHO}) then we have
\begin{eqnarray*}
u^{(p)}(t)=Q_{p}(t,u^{(k)}),\qquad t\geq t_0,
\end{eqnarray*}
where
\begin{eqnarray*}
Q_{p}(t,x)=\sum_{j=1}^{k-p}a_{j}\frac{(t-t_0)^{k-p-j}}{(k-p-j)!}+\int_{t_0}^{t}\frac{(t-s)^{k-p-1}}{(k-p-1)!}x(s)ds,\quad x\in
C([t_0,+\infty),{\R}),
\end{eqnarray*}
for every $p\in\overline{0,k-1}$. When $k=0$ we take $Q_{p}(t,x)\equiv x$.

So, for $x=u^{(k)}$, the equation (\ref{ecHO}) can be written as the integro-differential equation
\begin{eqnarray}
x^{(r+1)}=P(t,Q_{0}(t,x),Q_{1}(t,x),\dots,Q_{k-1}(t,x),x,x^{\prime},\dots,x^{(r)}).\label{HO1}
\end{eqnarray}

Further, by introducing the function $y=\mbox{col }(y_1,\cdots,y_{r+1})$ through the formulas
\begin{eqnarray*}
y_1=x,\thinspace y_2=x^{\prime},\dots, y_{r+1}=x^{(r)},
\end{eqnarray*}
the equation (\ref{HO1}) becomes (see \cite{Mustafa_3})
\begin{eqnarray*}
y^{\prime}=Ay+P[t,x,y],\qquad t\geq t_0,
\end{eqnarray*}
where
\begin{eqnarray*}
A=\left(
\begin{array}{lllll}
0 & 1 & 0 & \cdots & 0 \\
0 & 0 & 1 & \cdots & 0 \\
\vdots & \vdots & \vdots & \ddots & \vdots \\
0 & 0 & 0 & \cdots & 1 \\
0 & 0 & 0 & \cdots & 0
\end{array}
\right)\qquad\mbox{and}\qquad P[t,x,y]=\left(
\begin{array}{c}
0 \\
\vdots\\
0\\
P
\end{array}\right).
\end{eqnarray*}

The change of variables suggested in \cite{Mustafa_3}, namely $y=e^{At}z$, where $z=\mbox{col }(z_1,\dots,z_{r+1})$, leads to
\begin{eqnarray}
\left\{\begin{array}{ll}
y_{i}(t)=\sum\limits_{l=i}^{r+1}\frac{t^{l-i}}{(l-i)!}z_{l}(t),\\
z_{i}(t)=\sum\limits_{l=i}^{r+1}(-1)^{l-i}\frac{t^{l-i}}{(l-i)!}y_{l}(t)
\end{array}
\right.\label{HO2}
\end{eqnarray}
for $1\leq i\leq r+1$. We took into account that $(e^{At})^{-1}=e^{-At}$ when $t\geq t_0$.

We introduce the function $R(t,z)$ via the formula
\begin{eqnarray*}
R(t,z)=P\left(t,\left(Q_{p}\left(t,\sum_{q=1}^{r+1}\frac{t^{q-1}}{(q-1)!}z_{q}\right)\right)_{p\in\overline{0,k-1}},
\left(\sum_{l=d}^{r+1}\frac{t^{l-d}}{(l-d)!}z_{l}\right)_{d\in\overline{1,r+1}}\right).
\end{eqnarray*}

At this point, the relation $e^{At}z^{\prime}=P[t,x,y]$, cf. \cite{Mustafa_3}, becomes the integro-differential system
\begin{eqnarray}
\frac{dz_i}{dt}=(-1)^{r+1-i}\frac{t^{r+1-i}}{(r+1-i)!}R(t,z),\qquad 1\leq i\leq r+1.\label{sistintegrodifHO}
\end{eqnarray}

To establish the existence of a solution to the problem (\ref{ecHO}), (\ref{dateHO}) we impose a series of \textit{technical conditions} on both
the data $a_{k-i},c$ and the nonlinearity $P$.

Fix the numbers $M>0$, $s_{0}\in(0,1)$ and $s_{d}\in[s_0,1)$ such that $\sum\limits_{d=1}^{r+1}s_d=1$. Let also
$K:[t_0,+\infty)\rightarrow[0,+\infty)$ be a continuous function such that
\begin{eqnarray}
\int_{t_{0}}^{+\infty}K(s)ds\leq s_{0}M.\label{HO_hyp1}
\end{eqnarray}

We denote with $Z_{M}$ the subset of the space $[X_{2}(t_0;1)]^{r+1}$ (endowed with the product topology and respectively the product norm), see
\cite{Mustafa_3}, which consists of all functions $z(t)$ with the property that
\begin{eqnarray*}
\left\{
\begin{array}{ll}
\vert z_{1}(t)-c\vert\leq\frac{s_1}{s_0}\int_{t}^{+\infty}K(\tau)d\tau\\
\vert z_{i}(t)\vert\leq\frac{s_{i}}{s_0}\int_{t}^{+\infty}\frac{K(\tau)}{\tau^{i-1}}d\tau,\qquad i\in\overline{2,r+1},
\end{array}
\right.
\end{eqnarray*}
throughout $[t_0,+\infty)$.

The key technical condition of the Hale-Onuchic method is given by the inequality
\begin{eqnarray}
t^{r}\vert R(t,z)\vert\leq K(t),\qquad z\in Z_{M},\label{HO_hyp2}
\end{eqnarray}
for all $t\geq t_{0}$.

\begin{teorema}\label{teor1HO}\emph{(Hale, Onuchic, 1963, cf.
\cite[Theorem 1]{HaleOnuchic})} Assume that (\ref{HO_hyp1}), (\ref{HO_hyp2}) hold. Then, the system (\ref{sistintegrodifHO}) has at least one
solution in $Z_{M}$.
\end{teorema}

\textbf{Proof.} Consider the operator $T:Z_{M}\rightarrow[X_{2}(t_0;1)]^{r+1}$ with the formula $T(z)=w$, where $w=\mbox{col
}(w_1,\dots,w_{r+1})$ and
\begin{eqnarray*}
\left\{\begin{array}{ll}
w_{1}(t)=c-\int_{t}^{+\infty}(-1)^{r}\frac{s^{r}}{r!}R(s,z)ds\\
w_{i}(t)=-\int_{t}^{+\infty}(-1)^{r+1-i}\frac{s^{r+1-i}}{(r+1-i)!}R(s,z)ds,\qquad i\in\overline{2,r+1}.
\end{array}
\right.
\end{eqnarray*}

The set $Z_M$ is bounded, closed and convex. By adapting the technique from \cite{Mustafa_3}, we deduce that the operator $T$ is uniformly
continuous and the set $T(Z_M)$ is relatively compact.

It is easy to see that
\begin{eqnarray*}
\vert w_{1}(t)-c\vert\leq\int_{t}^{+\infty}K(\tau)d\tau\leq\frac{s_1}{s_0}\int_{t}^{+\infty}K(\tau)d\tau,
\end{eqnarray*}
and respectively
\begin{eqnarray*}
\vert w_{i}(t)\vert\leq\int_{t}^{+\infty}\frac{K(\tau)}{\tau^{i-1}}d\tau\leq\frac{s_i}{s_0}\int_{t}^{+\infty}\frac{K(\tau)}{\tau^{i-1}}d\tau,
\end{eqnarray*}
which implies $T(Z_{M})\subseteq Z_M$. The conclusion follows by application of the Schauder-Tikhonov fixed point theorem. $\square$

\begin{propozitia}\label{teor2HO}
Let $z\in Z_M$ be the solution of the system (\ref{sistintegrodifHO}) established in Theorem \ref{teor1HO}. Then, the functions
$(y_i)_{i\in\overline{1,r+1}}$ given by the first of formulas (\ref{HO2}) satisfy the estimates
\begin{eqnarray}
\vert y_{1}(t)-c\vert\leq\psi(t)\qquad\mbox{and}\qquad\vert y_{i+1}(t)\vert\leq\frac{\psi(t)}{t^i},\label{HO3}
\end{eqnarray}
where $\psi(t)=\frac{1}{s_0}\int_{t}^{+\infty}K(\tau)d\tau$, throughout $[t_0,+\infty)$. Furthermore,
\begin{eqnarray*}
y_{1}^{(j)}=y_{j+1},\qquad 0\leq j\leq r.
\end{eqnarray*}
\end{propozitia}

\textbf{Proof.} We have
\begin{eqnarray*}
\vert y_1(t)-c\vert\leq\vert z_1(t)-c\vert+\sum_{l=2}^{r+1}\frac{t^{l-1}}{(l-1)!}\vert
z_l(t)\vert\leq\psi(t)\sum_{l=1}^{r+1}\frac{s_l}{(l-1)!}\leq\psi(t),
\end{eqnarray*}
and respectively
\begin{eqnarray*}
\vert y_i(t)\vert\leq\frac{\psi(t)}{t^{i-1}}\sum_{l=i}^{r+1}\frac{t^{l-1}}{(l-1)!}\vert
z_l(t)\vert\leq\frac{\psi(t)}{t^{i-1}}\sum_{l=i}^{r+1}\frac{s_l}{(l-1)!}\leq\frac{\psi(t)}{t^{i-1}},\quad 2\leq i\leq r+1.
\end{eqnarray*}

Further, for $j\in\overline{1,r}$ fixed, we check that
\begin{eqnarray*}
y_{j}^{\prime}(t)&=&\sum_{q=j}^{r}\frac{t^{q-j}}{(q-j)!}z_{q+1}(t)+\sum_{l=j}^{r+1}\frac{t^{l-j}}{(l-j)!}z_{l}^{\prime}(t)\\
&=&y_{j+1}(t)+t^{r+1-j}R(t,z)\sum_{l=j}^{r+1}\frac{(-1)^{r+1+l}}{(l-j)!(r+1-l)!}=y_{j+1}(t),
\end{eqnarray*}
since
\begin{eqnarray*}
\sum_{l=j}^{r+1}\frac{(-1)^{r+1+l}}{(l-j)!(r+1-l)!}=\frac{1}{(r+1-j)!}\sum_{d=0}^{r+1-j}
\left(\begin{array}{c}d\\r+1-i\end{array}\right)(-1)^{d}=0
\end{eqnarray*}
according to Newton's binomial formula. $\square$

The conclusions of Proposition \ref{teor2HO} allow us to say that the problem (\ref{ecHO}), (\ref{dateHO}) has a solution $x(t)$. Moreover, the
relations (\ref{HO3}) describe the decay rate of the derivatives of this solution to  $c$, respectively to $0$ when $t\rightarrow+\infty$.

The condition (\ref{HO_hyp2}) is, in the general case, extremely complicated to verify. Certain simpler restrictions can be designed, fortunately, that
will lead to (\ref{HO_hyp2}).

\begin{teorema}\label{teor3HO}\emph{(\cite[Theorem 2]{HaleOnuchic})}
Consider the quasi-linear differential equation
\begin{eqnarray}
u^{(k+r+1)}+\sum_{j=1}^{k+r+1}f_{j}(t,u,u^{\prime},\dots,u^{(k+r)})u^{(k+r+1-j)}=0,\qquad t\geq1,\label{ecHO2}
\end{eqnarray}
where the functional coefficients $f_j:[1,+\infty)\times{\R}^{k+r+1}\rightarrow{\R}$ are continuous and verify the inequality
\begin{eqnarray*}
\vert f_{j}(t,u,u^{\prime},\dots,u^{(k+r)})\vert\leq h_{j}(t)L_{j}(\vert u\vert,\vert u^{\prime}\vert,\dots,\vert u^{(k+r)}\vert)
\end{eqnarray*}
throughout their entire domain of existence. Assume, further, that the comparison functions $L_{j}$ are continuous, monotone nondecreasing in
every argument and subjected to the restriction
\begin{eqnarray*}
\int^{+\infty}t^{j-1}h_{j}(t)L_{j}(\beta t^{k},\beta t^{k-1},\dots,\beta t^{-r})dt<+\infty
\end{eqnarray*}
for every $\beta>0$.

Then, given the numbers $c$, $(a_j)_{j\in\overline{1,k}}$, there exist $t_0\geq1$ and a solution $u(t)$ of equation (\ref{ecHO2}) defined in
$[t_0,+\infty)$ which satisfies the boundary conditions (\ref{dateHO}).
\end{teorema}

Let us mention at this point that results similar to Theorems \ref{teor1HO}, \ref{teor3HO} can be established also by following the fundamental technique due to Hartman and Onuchic \cite{HartmanOnuchic}. This method relies on a rather complicated (in the general case) change of variables and the admissibility theory devised by Massera, Sch\"{a}ffer and Corduneanu, see the discussion in \cite{Agarwal_et_al}. In fact, in \cite[Theorem 4]{Agarwal_et_al} we use this approach to deal with the existence of solutions that possess the same behavior as the one we look for in the present paper. 

\textbf{Proof.} (of Theorem \ref{teor3HO}) Fix $M>0$ and introduce $\beta=\vert c\vert+M+\sum\limits_{j=1}^{k}\vert a_{i}\vert$. Without providing a lower bound for $t_0$
at this time, we deduce that ($z\in Z_M$)
\begin{eqnarray}
&&\left\vert Q_{p}\left(t,\sum_{q=1}^{r+1}\frac{t^{q-1}}{(q-1)!}z_{q}(t)\right)\right\vert\leq\sum_{j=1}^{k-p}\vert a_j\vert
t^{k-p-1}+\int_{t_0}^{t}
\frac{(t-s)^{k-p-1}}{(k-p-1)!}\nonumber\\
&&\times\left(\sum_{q=1}^{r+1}\frac{s^{q-1}}{(q-1)1}\vert z_{q}(s)\vert\right)ds\leq\sum_{j=1}^{k-p}\vert a_j\vert t^{k-p-1}+\left[\vert
c\vert+M\left(\sum_{q=1}^{r+1}\frac{s_q}{(q-1)!}\right)\right]\nonumber\\
&&\times\int_{t_0}^{t} \frac{(t-s)^{k-p-1}}{(k-p-1)!}ds\leq\beta t^{k-p},\label{HO4}
\end{eqnarray}
and respectively
\begin{eqnarray}
&&\left\vert\sum_{l=d}^{r+1}\frac{t^{l-d}}{(l-d)!}z_{l}(t)\right\vert\leq\vert z_1(t)\vert+\sum_{l=d}^{r+1}\mbox{sgn
}(l-1)\frac{t^{1-d}}{(l-d)!}[t^{l-1}\vert z_l(t)\vert]\nonumber\\
&&\leq(\vert c\vert+M)t^{1-d}\leq\beta t^{1-d}\label{HO5}
\end{eqnarray}
for all $t\geq t_0$ and $0\leq p\leq k-1$, $1\leq d\leq r+1$.

Thus,
\begin{eqnarray}
&&t^{r}\vert R(t,z)\vert\nonumber\\
&&\leq\sum_{j=1}^{r+1}t^{r}\left\vert f_{j}\left(t,(Q_{p}(t,x))_{p\in\overline{0,k-1}},(x^{(i)})_{i\in\overline{0,r}}\right)\right\vert\vert
x^{(r+1-j)}(t)\vert\label{HO_suma1}\\
&&+\sum_{j=r+2}^{r+k+1}t^{r}\left\vert f_{j}\left(t,(Q_{p}(t,x))_{p\in\overline{0,k-1}},(x^{(i)})_{i\in\overline{0,r}}\right)\right\vert\vert
Q_{k+r+1-j}(t,x)\vert\label{HO_suma2},
\end{eqnarray}
where, according to Proposition \ref{teor2HO}, we have
\begin{eqnarray*}
x^{(i)}(t)=y_{i+1}(t)=\sum_{l=i+1}^{r+1}\frac{t^{l-i-1}}{(l-i-1)!}z_{l}(t).
\end{eqnarray*}

The estimates (\ref{HO4}), with $p=k+r+1-j$, and (\ref{HO5}), with $d=r+2-j$, allow us to majorize any of the terms from the sum
(\ref{HO_suma1}), and respectively of the sum (\ref{HO_suma2}) with the quantity (see \cite[p. 68]{HaleOnuchic})
\begin{eqnarray*}
\beta t^{j-1}h_{j}(t)L_{j}(\beta t^{k},\beta t^{k-1},\dots,\beta t^{-r}),\qquad t\geq t_0.
\end{eqnarray*}

The function $K:[t_0,+\infty)\rightarrow[0,+\infty)$ is now the sum
\begin{eqnarray*}
K(t)=\beta\sum_{j=1}^{k+r+1}t^{j-1}h_{j}(t)L_{j}(\beta t^{k},\beta t^{k-1},\dots,\beta t^{-r}),
\end{eqnarray*}
where $t_0\geq1$ verifies the restriction (\ref{HO_hyp1}). $\square$

The significance of (\ref{HO_hyp2}) is that, in a series of cases which cannot be discussed by means of Theorem \ref{teor3HO}, \textit{the
nonlinearity of the differential equation (\ref{ecHO}) doesn't have to be necessarily "small" but it has only to become "small" along some of
the solutions of the equation}, cf. \cite[p. 61]{HaleOnuchic}. See also the discussion in \cite[p. 363]{Agarwal_et_al}.

General results about the existence throughout unbounded intervals of the solutions to various boundary value problems associated to nonlinear
differential equations can be found in \cite{Kartsatos1974,Kartsatoscarte_2}. We mention, in this respect, the monographs
\cite{Brauer_carte,Coppel_integr} and the papers \cite{Cecchi_pctfix,Talpalaru}.

Similar investigations, regarding the issue of disconjugacy, for quasi-linear differential equations can be read in
\cite{Schuur1983,Schuur1979,Schuur1982}. The existence of special families of solutions (Chebyshev or Descartes systems) is established via
fixed point theories of multi-functions (the Fan-Glicksberg theorem, \cite[p. 82]{Schuur1979}). This type of approach in the asymptotic
integration of quasi-linear ordinary differential equations is attributed to Kartsatos, cf. \cite[p. 911]{Schuur1982}.

Some interesting details regarding the existence over bounded intervals of solutions to general boundary value problems can be read in the
contribution by Avramescu \cite{Avramescu_Annaligen}.

\section{Local theories for the asymptotic integration of equation (\ref{ecEF1})}

Let us come back to the Emden-Fowler equation (\ref{ecEF1}). We shall analyze in the following the asymptotic features of some of its
(non-oscillatory) solutions using Hale-Onuchic results.

By \textit{local theories} we mean, in the spirit of condition (\ref{HO_hyp2}), those results in whose hypotheses \textit{the behavior of the
nonlinearity of a differential equation is described on a given family of functions} (with these functions we compare either the solution we are
looking for or other functional quantities associated with the solution) and not throughout its entire domain of existence.

The results are based on an observation made by Dub\'{e} and Mingarelli \cite[Eq. (2.1)]{DubeMingarelli_HaleOnuchic}, according to which
\textit{one can combine the hypothesis (\ref{HO_hyp2}) with a (local) Lipschitz restriction upon the nonlinearity of the differential equation}.
In this way, the functional analysis involved in the investigation is greatly simplified.

\begin{teorema}\label{teorDubeMingarelli_EF}\emph{(Dub\'{e}, Mingarelli, 2004, cf. \cite[Theorem
2.1]{DubeMingarelli_HaleOnuchic})} Let $f:[t_0,+\infty)\times{\R}\rightarrow[0,+\infty)$ be a continuous function such that
\begin{eqnarray}
\int_{t_0}^{+\infty}(t-t_0)f(t,u(t))dt\leq M\label{Dube1}
\end{eqnarray}
for all $u\in X_{M}$, where $M>0$ is fixed and
\begin{eqnarray*}
X_{M}=\{u\in C([t_0,+\infty),{\R}):0\leq u(t)\leq M\mbox{ for all }t\geq t_0\}.
\end{eqnarray*}

Assume that there exists the continuous function $k:[t_0,+\infty)\rightarrow[0,+\infty)$ subjected to
\begin{eqnarray}
\eta=\int_{t_0}^{+\infty}(t-t_0)k(t)dt<1\label{Dube_hyp1}
\end{eqnarray}
and
\begin{eqnarray}
\vert f(t,u_{2}(t))-f(t,u_{1}(t))\vert\leq k(t)\vert u_{2}(t)-u_{1}(t)\vert,\qquad t\geq t_0,\thinspace u_{1,2}\in X_{M}.\label{Dube_hyp2}
\end{eqnarray}

Then, the differential equation
\begin{eqnarray}
x^{\prime\prime}+f(t,x)=0,\qquad t\geq t_0\geq0,\label{ec_DubeMingarelli}
\end{eqnarray}
has a solution $x(t)$ defined in $[t_0,+\infty)$ with the asymptotic profile given by
\begin{eqnarray}
\lim\limits_{t\rightarrow+\infty}x(t)=M.\label{alura_DubeMingarelli}
\end{eqnarray}
\end{teorema}

\textbf{Proof.} We introduce the integral operator $T:X_{M}\rightarrow C([t_0,+\infty),{\R})$ with the formula
\begin{eqnarray*}
T(u)(t)=M-\int_{t}^{+\infty}(s-t)f(s,u(s))ds,\qquad u\in X_{M},
\end{eqnarray*}
for all $t\geq t_0$. It is easy to deduce that, by a double differentiation with respect to $t$, any  (eventual) fixed point of operator $T$
verifies the equation (\ref{ec_DubeMingarelli}) throughout $[t_0,+\infty)$.

Obviously, $T(X_{M})\subseteq X_{M}$ according to the Hale-Onuchic condition (\ref{Dube1}).

Given the pairs $(u_1,u_2)$ of elements from $X_M$, we introduce the natural distance $d$ by
\begin{eqnarray}
d(u_1,u_2)=\sup\limits_{t\geq t_0}\{\vert u_1(t)-u_2(t)\vert\}.\label{distDube}
\end{eqnarray}
It can be established in a standard manner, see \cite[Chapter I]{Kartsatoscarte_2}, that the metric space $S_{M}=(X_M,d)$ is complete.

The next estimates
\begin{eqnarray}
\vert T(u_2)(t)-T(u_1)(t)\vert&\leq&\int_{t}^{+\infty}(s-t)\vert
f(s,u_2(s))-f(s,u_1(s))\vert ds\nonumber\\
&\leq&\int_{t}^{+\infty}(s-t_0)k(s)\vert u_2(s)-u_1(s)\vert ds\nonumber\\
&\leq&\int_{t}^{+\infty}(s-t_0)k(s)ds\cdot d(u_2,u_1),\label{estimDubeoperator}
\end{eqnarray}
where $u_{1,2}\in X_M$ \c{s}i $t\geq t_0$, implies $d(T(u_2),T(u_1))\leq\eta d(u_1,u_2)$, that is \textit{the operator $T$ is a contraction in
$S_M$}.

The Banach contraction principle \cite[pp. 19--20]{Kartsatoscarte_2} ensures that the equation (\ref{ec_DubeMingarelli}) has the solution we
were looking for. $\square$

In its essence, the Hale-Onuchic "philosophy" of asymptotic integration of nonlinear differential equations reduces to \textit{transforming the
boundary value problem into the existence problem of a fixed point to an integral operator $T$}, a common fact in this field, and to
\textit{identifying an invariant set ($Z_M$, $X_M$, etc)} in which one can use, with minimal effort, a fixed point theorem. This is why, in a
Hale-Onuchic type of approach, the verification of the hypotheses of such a theorem is quite easy, the weight leaning upon the formula of the
integral operator (associated to an intermediate integro-differential problem), respectively upon the detection of invariant sets. Some authors,
maybe too drastic in this respect, exclude the verification of the hypotheses to the fixed point theorem from the investigation, e.g. \cite[p.
v]{Eastham_integr}.

By taking $M=1$, $k(t)\equiv\lambda q(t)$ and $f(t,u)\equiv q(t)u^{\lambda}$, Theorem \ref{teorDubeMingarelli_EF} establishes the existence of a
solution to the boundary value problem (\ref{bvp_Atkinson}). Atkinson's original claim is improved in

\begin{teorema}\label{teorAtkinson_imbunatatita}
Assume that (\ref{Atkinson_hyp1}) holds. Then, there exists $p>1$ such that, for all $c\in(0,1]$, the equation (\ref{ecEF1}) has at least one
solution $x(t)$ defined in $[t_0,+\infty)$ with the property that
\begin{eqnarray}
\frac{c}{p}\leq x(t)\leq c,\qquad t\geq t_0.\label{estim_Atkinsonimbunatatita}
\end{eqnarray}

The asymptotic profile of the solution is given by $x(t)=c+o(1)$, respectively $x^{\prime}(t)=o(t^{-1})$ when $t\rightarrow+\infty$.
\end{teorema}

\textbf{Proof.} Introduce the complete metric space $S=(X,d)$, where
\begin{eqnarray*}
X=\left\{u\in C([t_0,+\infty),{\R}):\frac{c}{p}\leq u(t)\leq c\mbox{ for all }t\geq t_0\right\}
\end{eqnarray*}
and the distance $d$ is given by (\ref{distDube}). Here, the number $p>1$ is fixed such that to have ($\eta\in(0,1)$!)
\begin{eqnarray}
\int_{t_0}^{+\infty}tq(t)dt<\frac{1}{\lambda}\leq\frac{p-1}{p}<1.\label{Mustafa_estim_p}
\end{eqnarray}

The nonlinearity $f(t,u)\equiv q(t)u^{\lambda}$ satisfies the condition (\ref{Dube_hyp2}) for $k(t)\equiv\lambda q(t)$. Also, the operator
$T:X\rightarrow C([t_0,+\infty),{\R})$ with the formula
\begin{eqnarray}
T(u)(t)=c-\int_{t}^{+\infty}(s-t)q(s)[u(s)]^{\lambda}ds,\qquad u\in X,\thinspace t\geq t_{0},\label{primaformopameliorat}
\end{eqnarray}
verifies the estimate (\ref{estimDubeoperator}), meaning that it is a contraction of coefficient $\eta$.

We have
\begin{eqnarray*}
c&\geq&
c-\int_{t}^{+\infty}(s-t)q(s)\left(\frac{c}{p}\right)^{\lambda}ds\\
&\geq&
T(u)(t)=c-\int_{t}^{+\infty}(s-t)q(s)[u(s)]^{\lambda}ds\\
&\geq&c-c^{\lambda}\int_{t}^{+\infty}(s-t)q(s)ds\geq
c\left[1-\int_{t_0}^{+\infty}sq(s)ds\right]\\
&\geq&\frac{c}{p},\qquad u\in X,\thinspace t\geq t_0,
\end{eqnarray*}
that is $T(X)\subseteq X$.

According to the contraction principle, the operator $T$  has a fixed point in $X$, denoted $x(t)$. Obviously, $x(t)=c+o(1)$ when
$t\rightarrow+\infty$. Further,
\begin{eqnarray*}
0\leq x^{\prime}(t)=\int_{t}^{+\infty}q(s)[x(s)]^{\lambda}ds\leq c^{\lambda}\int_{t}^{+\infty}q(s)ds\leq\int_{t}^{+\infty}q(s)ds=o(t^{-1})
\end{eqnarray*}
as $t\rightarrow+\infty$. See \cite{MustafaNA2005}. $\square$

\begin{observatia}
\emph{We impose this restriction upon $c$ just to make use of (\ref{Atkinson_hyp1}). Elsewhere, in the spirit of \cite{CoffWong}, fix $C$ such
that $C>\lambda(C-\vert c\vert)>0$ and $C^{\lambda}\int_{t_0}^{+\infty}t\vert q(t)\vert dt\leq C-\vert c\vert$. Then, the operator
$T:{\cal{C}}\rightarrow X_{2}(t_0;1)$, see \cite{Mustafa_3}, with the formula given by (\ref{primaformopameliorat}) and ${\cal{C}}=\{u\in
X_{2}(t_0;1):\vert u(t)\vert\leq C\mbox{ for all }t\geq t_0\}$, is a contraction of coefficient $\lambda\frac{C-\vert c\vert}{C}$. Its fixed
point in ${\cal{C}}$ is the solution we are looking for.}
\end{observatia}

The conclusion of Theorem \ref{teorAtkinson_imbunatatita} can be reached as well by a different local theory. To this end, we establish a
theorem that improves upon the results in \cite[Theorem 1]{MustafaNA2005}.

\begin{teorema}\label{teormustafa1}
Set $M\in{\R}$ and let $\alpha,\beta:[t_0,+\infty)\rightarrow{\R}$ be continuous functions, absolutely integrable over $[t_0,+\infty)$, with
$\alpha(t)\leq\beta(t)$ for all $t\geq t_0$ and such that
$\lim\limits_{t^{\prime}\rightarrow+\infty}\alpha(t^{\prime})=\lim\limits_{t^{\prime}\rightarrow+\infty}\beta(t^{\prime})=0$.

Given the sets
\begin{eqnarray*}
C_M&=&\left\{u\in C([t_0,+\infty),{\R}):M-\int_{t}^{+\infty}\beta(s)ds\leq
u(t)\right.\\
&&\left.\leq M-\int_{t}^{+\infty}\alpha(s)ds\mbox{ for all }t\geq t_0\right\}
\end{eqnarray*}
and
\begin{eqnarray*}
D=\{v\in C([t_0,+\infty),{\R}):\alpha(t)\leq v(t)\leq\beta(t)\mbox{ for all }t\geq t_0\},
\end{eqnarray*}
assume that the double inequality takes place
\begin{eqnarray}
\alpha(t)\leq\int_{t}^{+\infty}f(s,u(s),v(s))ds\leq\beta(t),\qquad t\geq t_0,\label{condMustafa1}
\end{eqnarray}
for all $u\in C_M$ and $v\in D$, where the function $f:[t_0,+\infty)\times{\R}^{2}\rightarrow{\R}$ is continuous. As a plus,
\begin{eqnarray*}
&&\vert f(t,u_{2}(t),v_{2}(t))-f(t,u_{1}(t),v_{1}(t))\vert\\
&&\leq k_1(t)\vert u_{2}(t)-u_{1}(t)\vert+k_2(t)\vert v_{2}(t)-v_1(t)\vert
\end{eqnarray*}
for all $u_{1,2}\in C_M$, $v_{1,2}\in D$ and $t\geq t_0$. Here, the functions $k_{1,2}:[t_0,+\infty)\rightarrow[0,+\infty)$ are continuous and
there exists a number $\zeta>0$ such that
\begin{eqnarray*}
&&\chi=\zeta\int_{t_0}^{+\infty}k_{1}(t)dt+\int_{t_0}^{+\infty}(t-t_0)k_{1}(t)dt\\
&&+\int_{t_0}^{+\infty}k_{2}(t)dt+\frac{1}{\zeta}\int_{t_0}^{+\infty}(t-t_0)k_2(t)dt<1.
\end{eqnarray*}

Then, the boundary value problem
\begin{eqnarray}
\left\{
\begin{array}{ll}
x^{\prime\prime}+f(t,x,x^{\prime})=0,\qquad t\geq t_0\geq0,\\
\lim\limits_{t\rightarrow+\infty}x(t)=M\\
\alpha(t)\leq x^{\prime}(t)\leq\beta(t),\qquad t\geq t_0,
\end{array}
\right.\label{bvp_Mustafa1}
\end{eqnarray}
has a unique solution.
\end{teorema}

\textbf{Proof.} Define the distance $d$ between the elements $v_1$ and $v_2$ of the set $D$ by the formula
\begin{eqnarray*}
d(v_1,v_2)=\Vert v_1-v_2\Vert_{L^{1}((t_0,+\infty),{\R})}+\zeta\sup\limits_{t\geq t_{0}}\{\vert v_1(t)-v_2(t)\vert\}.
\end{eqnarray*}
The dominated convergence theorem, see \cite[pp. 20--21]{Evans_carte}, ensures that the metric space $S=(D,d)$ is complete.

We introduce the operator $T:D\rightarrow C([t_0,+\infty),{\R})$ via the formula
\begin{eqnarray*}
T(v)(t)=\int_{t}^{+\infty}f\left(s,M-\int_{s}^{+\infty}v(\tau)d\tau,v(s)\right)ds,\qquad v\in D,\thinspace t\geq t_0.
\end{eqnarray*}
The restriction (\ref{condMustafa1}) shows that $T(D)\subseteq D$ since $M-\int_{(\cdot)}^{+\infty}v(s)ds\in C_{M}$ for all $v\in D$.

\begin{afirmatia}
The operator $T:D\rightarrow D$ is a contraction with coefficient $\chi$.
\end{afirmatia}

In fact, we have
\begin{eqnarray*}
\zeta\vert T(v_2)(t)-T(v_1)(t)\vert&\leq&\zeta\int_{t}^{+\infty}k_1(s)\int_{s}^{+\infty}\vert
v_2(\tau)-v_1(\tau)\vert d\tau ds\\
&+&\int_{t}^{+\infty}k_2(s)[\zeta\vert v_{2}(s)-v_1(s)\vert]ds\\
&\leq&\left[\zeta\int_{t_0}^{+\infty}k_1(s)ds+\int_{t_0}^{+\infty}k_2(s)ds\right]d(v_1,v_2)
\end{eqnarray*}
and
\begin{eqnarray*}
&&\int_{t}^{+\infty}\vert T(v_2)(s)-T(v_1)(s)\vert
ds\leq\int_{t}^{+\infty}(s-t)k_1(s)\\
&&\times\int_{s}^{+\infty}\vert v_2(\tau)-v_1(\tau)\vert d\tau ds
+\frac{1}{\zeta}\int_{t}^{+\infty}(s-t)k_2(s)[\zeta\vert v_{2}(s)-v_1(s)\vert]ds\\
&&\leq\left[\int_{t_0}^{+\infty}(s-t_0)k_1(s)ds+\frac{1}{\zeta}\int_{t_0}^{+\infty}(s-t_0)k_2(s)ds\right]d(v_1,v_2).
\end{eqnarray*}

Now,
\begin{eqnarray*}
\int_{t}^{+\infty}\vert T(v_2)(s)-T(v_1)(s)\vert ds+\zeta\vert T(v_2)(t)-T(v_1)(t)\vert\leq\chi d(v_1,v_2)
\end{eqnarray*}
for all $v_{1,2}\in D$ and $t\geq t_0$, which validates the Claim.

Denote with $v_{0}(t)$ the fixed point of operator $T$ in $D$. Then, the function $x(t)\equiv M-\int_{t}^{+\infty}v_{0}(s)ds$ is the solution we
are looking for. $\square$

\begin{corolarul}Let $g:[t_0,+\infty)\rightarrow[0,+\infty)$
be a continuous function, integrable over $[t_0,+\infty)$ and such that $\lim\limits_{t\rightarrow+\infty}g(t)=0$. Assume as well that, in the
statement of Theorem \ref{teorDubeMingarelli_EF}, the hypothesis (\ref{Dube1}) is replaced by the restrictions
\begin{eqnarray*}
\int_{t}^{+\infty}f(s,u(s))ds\leq g(t),\thinspace u\in X_M,\thinspace t\geq
t_0,\qquad\mbox{and}\qquad\int_{t_0}^{+\infty}g(t^{\prime})dt^{\prime}\leq M.
\end{eqnarray*}

Then, if the function $f(t,x)$ verifies the conditions (\ref{Dube_hyp1}), (\ref{Dube_hyp2}), there will be a solution $x(t)$ of equation
(\ref{ec_DubeMingarelli}) with the asymptotic profile (\ref{alura_DubeMingarelli}) such that
\begin{eqnarray*}
0\leq x^{\prime}(t)\leq g(t),\qquad t\geq t_0.
\end{eqnarray*}
\end{corolarul}

\textbf{Proof.} Take $\alpha=0$, $\beta=g$, $k_1=k$, $k_2=0$ and $\zeta\in(0,1)$ with the property that
\begin{eqnarray}
\zeta<(1-\eta)\left(\int_{t_0}^{+\infty}k(t)dt\right)^{-1}.\label{estimzeta}
\end{eqnarray}

Since $\chi=\zeta\int_{t_0}^{+\infty}k(t)dt+\eta<1$, the conclusion follows readily from  Theorem \ref{teormustafa1}. $\square$

\vspace{0.5cm}

\textbf{A second proof of Theorem \ref{teorAtkinson_imbunatatita}.} We shall use Theorem \ref{teormustafa1}. So, fix the numbers $p>1$,
$\zeta\in(0,1)$ in order for (\ref{Mustafa_estim_p}), (\ref{estimzeta}) to hold.

Define the functions
\begin{eqnarray*}
\alpha(t)=\left(\frac{c}{p}\right)^{\lambda}\int_{t}^{+\infty}q(s)ds,\qquad\beta(t)=c^{\lambda}\int_{t}^{+\infty}q(s)ds,\qquad t\geq t_0.
\end{eqnarray*}
Also, introduce $k_1=k$, $k_2=0$.

We have
\begin{eqnarray}
\int_{t_0}^{+\infty}\beta(t)dt=c^{\lambda}\int_{t_0}^{+\infty}(t-t_0)q(t)dt\leq
c\int_{t_0}^{+\infty}tq(t)dt<c\left(1-\frac{1}{p}\right).\label{estimmustafa2}
\end{eqnarray}

According to \cite[p. 183]{MustafaNA2005}, and taking into account (\ref{estimmustafa2}), the double inequality (\ref{condMustafa1}) follows
from the next estimates
\begin{eqnarray*}
&&\alpha(t)=\int_{t}^{+\infty}q(s)\left(\frac{c}{p}\right)^{\lambda}ds\leq\int_{t}^{+\infty}q(s)
\left(c-\int_{s}^{+\infty}\beta(\tau)d\tau\right)^{\lambda}ds\\
&&\leq\int_{t}^{+\infty}q(s)[u(s)]^{\lambda}ds=\int_{t}^{+\infty}f(s,u(s))ds,\quad
u\in C_c,\thinspace (M=c\in(0,1]!)\\
&&\leq\int_{t}^{+\infty}q(s) \left(c-\int_{s}^{+\infty}\alpha(\tau)d\tau\right)^{\lambda}ds\leq c^{\lambda}\int_{t}^{+\infty}q(s)ds=\beta(t)
\end{eqnarray*}
for all $t\geq t_0$.

Consequently, by applying Theorem \ref{teormustafa1}, we establish the existence of a solution with the formula $x(t)\equiv
c-\int_{t}^{+\infty}v_{0}(s)ds$ of the boundary value problem (\ref{bvp_Mustafa1}), where $v_0$ is the fixed point of operator $T$ in $D$.

The estimate (\ref{estim_Atkinsonimbunatatita}) is a by-product of (\ref{estimmustafa2}). More precisely,
\begin{eqnarray*}
&&\frac{c}{p}\leq c\left[1-\int_{t_0}^{+\infty}(t-t_0)q(t)dt\right]\leq
c-c^{\lambda}\int_{t_0}^{+\infty}\int_{t^{\prime}}^{+\infty}q(s)dsdt^{\prime}\\
&&\leq c-\int_{t}^{+\infty}\beta(t^{\prime})dt^{\prime}\leq c-\int_{t}^{+\infty}v_{0}(t^{\prime})dt^{\prime}=x(t)\leq c,\qquad t\geq t_0.
\end{eqnarray*}

The proof is complete. $\square$

Another local result regards the \textit{linear-like} solutions of equation (\ref{ecEF1}), cf. \cite[p. 356]{Agarwal_et_al}.

\begin{teorema}\emph{(\cite[Theorem 2.1]{Mustafa_mexicana})} Take
$t_0\geq1$, $A,x_0\in{\R}$, $\nu\in[0,1)$ and the continuous functions $\alpha,\beta:[t_0,+\infty)\rightarrow{\R}$ such that
$\alpha(t^{\prime})\leq\beta(t^{\prime})$ for all $t^{\prime}\geq t_0$ and $\alpha(t),\beta(t)=o(t^{-\nu})$ when $t\rightarrow+\infty$.

Consider the set $E_{A,x_0}$ given by
\begin{eqnarray*}
E_{A,x_0}&=&\{u\in C^{1}([t_0,+\infty),{\R}):u(t_0)=x_0,\thinspace\alpha(t)\leq
u^{\prime}(t)-A\\
&&\leq\beta(t)\mbox{ for all }t\geq t_0\}
\end{eqnarray*}
and assume that
\begin{eqnarray*}
\alpha(t)\leq\int_{t}^{+\infty}f(s,u(s))ds\leq\beta(t),\qquad u\in E_{A,x_0},\thinspace t\geq t_0,
\end{eqnarray*}
where the continuous function $f:[t_0,+\infty)\times{\R}\rightarrow{\R}$ is the nonlinearity of the equation (\ref{ec_DubeMingarelli}). As a
plus,
\begin{eqnarray*}
\vert f(t,u_{2}(t))-f(t,u_{1}(t))\vert\leq\frac{k(t)}{t}\vert u_{2}(t)-u_{1}(t)\vert
\end{eqnarray*}
for all $u_{1,2}\in E_{A,x_0}$ and $t\geq t_0$. Here, the function $k:[t_0,+\infty)\rightarrow[0,+\infty)$ is continuous and such that
\begin{eqnarray}
\varpi=\frac{1}{1-\nu}\int_{t_0}^{+\infty}k(t)dt<1.\label{estimmustafa3}
\end{eqnarray}

Then, the boundary value problem
\begin{eqnarray*}
\left\{
\begin{array}{ll}
x^{\prime\prime}+f(t,x)=0,\qquad t\geq t_0,\\
x(t_0)=x_0,\\
x(t)=At+o(t^{1-\nu})\qquad\mbox{when }t\rightarrow+\infty
\end{array}
\right.
\end{eqnarray*}
has a unique solution $x(t)$ with the property that
\begin{eqnarray*}
\alpha(t)\leq x^{\prime}(t)-A\leq\beta(t),\qquad t\geq t_0.
\end{eqnarray*}

In particular, if $\int^{+\infty}\alpha(t)dt=+\infty$, then $\lim\limits_{t\rightarrow+\infty}[x(t)-At]=+\infty$.
\end{teorema}

\textbf{Proof.} Define the distance between the elements $u_1$ and $u_2$ of set $E_{A,x_0}$ by the formula
\begin{eqnarray*}
d(u_1,u_2)=\sup\limits_{t\geq t_0}\{t^{\nu}\vert u_{1}^{\prime}(t)-u_{2}^{\prime}(t)\vert\}.
\end{eqnarray*}
One can see easily that the metric space $S=(E_{A,x_0},d)$ is complete.

Let the operator $T:E_{A,x_0}\rightarrow C^{1}([t_0,+\infty),{\R})$ be given by
\begin{eqnarray*}
T(u)(t)=x_0+A(t-t_0)+\int_{t_0}^{t}\int_{s}^{+\infty}f(\tau,u(\tau))d\tau ds,\qquad u\in E_{A,x_0},\thinspace t\geq t_0.
\end{eqnarray*}
Evidently, $T(E_{A,x_0})\subseteq E_{A,x_0}$.

From the next estimates, namely
\begin{eqnarray*}
\vert [T(u_2)]^{\prime}(t)-[T(u_1)]^{\prime}(t)\vert&\leq&\int_{t}^{+\infty}\frac{k(s)}{s}\vert
u_{2}(s)-u_1(s)\vert ds\\
&\leq&\int_{t}^{+\infty}\frac{k(s)}{s}\int_{t_0}^{s}\vert
u_{2}^{\prime}(\tau)-u_{1}^{\prime}(\tau)\vert d\tau ds\\
&\leq&\int_{t}^{+\infty}\frac{k(s)}{s}\int_{t_0}^{s}\frac{d\tau}{{\tau}^{\nu}}ds\cdot
d(u_1,u_2)\\
&\leq& t^{-\nu}\left(\frac{1}{1-\nu}\int_{t_0}^{+\infty}k(s)ds\right)d(u_1,u_2),
\end{eqnarray*}
we conclude that the operator $T$ is a contraction of coefficient $\varpi$ (we recall the restriction (\ref{estimmustafa3})!). Its fixed point
in $E_{A,x_0}$, denoted $x$, is the solution we have searched for.

The relations
\begin{eqnarray*}
x(t)=T(x)(t)\geq x_{0}+A(t-t_0)+\int_{t_0}^{t}\alpha(s)ds,\qquad t\geq t_0,
\end{eqnarray*}
allow us to prove the last requirement of the theorem. $\square$

\begin{corolarul}\label{cormustafa1} Set $t_0,\lambda\geq1$, $\nu\in[0,1)$ and $A>0$. Assume that
the continuous function $q:[t_0,+\infty)\rightarrow[0,+\infty)$ verifies the conditions
\begin{eqnarray}
\lambda c_{\nu}(A+c_{\nu})^{\lambda-1}<1-\nu\quad\mbox{and}\quad
\int_{t_0}^{+\infty}t^{\lambda}q(t)dt<\frac{c_{\nu}}{(A+c_{\nu})^{\lambda}},\label{estimmustafa4}
\end{eqnarray}
where $c_{\nu}=\int_{t_0}^{+\infty}t^{\lambda+\nu}q(t)dt$.

Then, the equation (\ref{ecEF1}) possesses the solution $x(t)$ with the asymptotic profile $x(t)=At+\omega(t)=At+o(t^{1-\nu})$ as
$t\rightarrow+\infty$, where
\begin{eqnarray*}
A^{\lambda}\int_{t_0}^{t}\int_{s}^{+\infty}\tau^{\lambda}q(\tau)d\tau ds\leq\omega(t)\leq
(A+c_{\nu})^{\lambda}\int_{t_0}^{t}\int_{s}^{+\infty}\tau^{\lambda}q(\tau)d\tau ds
\end{eqnarray*}
for all $t\geq t_0$.
\end{corolarul}

A simplification of hypothesis (\ref{estimmustafa4}) can be found in (\cite[Corollary 2.4]{Mustafa_mexicana}). For another variant, assume that
$\left(2-\frac{1}{\lambda}\right)\nu<1$.

\begin{afirmatia}
The conditions
\begin{eqnarray*}
\lambda\int_{t_0}^{+\infty}t^{\lambda+\left(2-\frac{1}{\lambda}\right)\nu}q(t)dt<1-\nu\quad\mbox{and}\quad (A+1)^{\lambda}<t_{0}^{\nu}
\end{eqnarray*}
imply (\ref{estimmustafa4}).
\end{afirmatia}

In fact, as $\lambda,t_0\geq1$, we have
\begin{eqnarray*}
c_{\nu}\leq\lambda t_{0}^{\nu\frac{\lambda-1}{\lambda}}\int_{t_0}^{+\infty}t^{\lambda+\nu}q(t)dt\leq\lambda
\int_{t_0}^{+\infty}t^{\lambda+\left(2-\frac{1}{\lambda}\right)\nu}q(t)dt<1-\nu<1,
\end{eqnarray*}
and respectively
\begin{eqnarray*}
\lambda c_{\nu}(A+c_{\nu})^{\lambda-1}<\lambda c_{\nu}(A+1)^{\lambda-1}\leq\lambda
t_{0}^{\nu\frac{\lambda-1}{\lambda}}\int_{t_0}^{+\infty}t^{\lambda+\nu}q(t)dt<1-\nu.
\end{eqnarray*}

With regard to the second of restrictions (\ref{estimmustafa4}), notice that
\begin{eqnarray*}
\int_{t_{0}}^{+\infty}t^{\lambda}q(t)dt&\leq&
t_{0}^{-\nu}\int_{t_0}^{+\infty}t^{\lambda+\nu}q(t)dt=t_{0}^{-\nu}c_{\nu}<c_{\nu}(A+1)^{-\lambda}\\
&<&c_{\nu}(A+c_{\nu})^{-\lambda}.
\end{eqnarray*}

The Claim is validated.

\vspace{0.5cm}

\textbf{Proof of Corollary \ref{cormustafa1}.} Fix $x_0=At_0$. Introduce the functions
\begin{eqnarray*}
\alpha(t)=A^{\lambda}\int_{t}^{+\infty}s^{\lambda}q(s)ds,\qquad\beta(t)=(A+c_\nu)^{\lambda}\int_{t}^{+\infty}s^{\lambda}q(s)ds,
\end{eqnarray*}
where $t\geq t_0$. The second of relations (\ref{estimmustafa4}) implies $\beta(t)\leq c_\nu$ in $[t_0,+\infty)$. Also, it is obvious that we
have  $\alpha(t),\beta(t)=o(t^{-\nu})$ when $t\rightarrow+\infty$. This is a consequence of the estimate
$\beta(t)\leq(A+c_\nu)^{\lambda}t^{-\nu}\int_{t}^{+\infty}s^{\lambda+\nu}q(s)ds$.

Given $u\in E_{A,x_0}$, we have
\begin{eqnarray*}
&&\alpha(t)\leq\int_{t}^{+\infty}q(s)\left(As+\int_{t_0}^{s}\alpha(\tau)d\tau\right)^{\lambda}ds\leq\int_{t}^{+\infty}q(s)[u(s)]^{\lambda}ds\\
&&=[T(u)]^{\prime}(t)-A\leq\int_{t}^{+\infty}q(s)\left(As+\int_{t_0}^{s}\beta(\tau)d\tau\right)^{\lambda}ds\\
&&\leq\int_{t}^{+\infty}q(s)[(A+c_{\nu})s]^{\lambda}ds=\beta(t),\qquad t\geq t_0.
\end{eqnarray*}

Further, via the mean value theorem, we get
\begin{eqnarray*}
&&\vert f(t,u_2(t))-f(t,u_1(t))\vert=t^{\lambda}q(t)\left\vert\left(\frac{u_2(t)}{t}\right)^{\lambda}-
\left(\frac{u_1(t)}{t}\right)^{\lambda}\right\vert\\
&&\leq\lambda t^{\lambda}q(t)\left[\frac{1}{t}\left(At+\int_{t_0}^{t}\beta(s)ds\right)\right]^{\lambda-1}\frac{\vert
u_2(t)-u_1(t)\vert}{t}\\
&&\leq\frac{\lambda(A+c_\nu)^{\lambda-1}t^{\lambda}q(t)}{t}\vert u_2(t)-u_1(t)\vert,\qquad u_{1,2}\in E_{A,x_0},
\end{eqnarray*}
throughout $[t_0,+\infty)$.

The first relation in (\ref{estimmustafa4}) implies (\ref{estimmustafa3}) for $k(t)\equiv\lambda(A+c_\nu)^{\lambda-1}t^{\lambda}q(t)$, by taking
into account  that $t_0\geq1$. $\square$

The following local theory is devoted to the condition (\ref{wronsk_sign}).

\begin{teorema}\label{MustRogoPDES}\emph{(\cite[Theorem 2.3]{MustafaRogovPDEs})} Consider
$t_0\geq1$, $a,b\geq0$, $c\in(0,1]$ and the bounded continuous functions $\alpha,\beta:[t_0,+\infty)\rightarrow[0,+\infty)$ such that
$\alpha(t^{\prime})\leq\beta(t^{\prime})$ for all $t^{\prime}\geq t_0$.

Introduce the set $F_{a,b,c}$ by the formula
\begin{eqnarray*}
F_{a,b,c}&=&\left\{u\in
C([t_0,+\infty),{\R}):at+b+t\int_{t}^{+\infty}\frac{\alpha(s)}{s^{1+c}}ds\leq u(t)\right.\\
&&\left.\leq at+b+t\int_{t}^{+\infty}\frac{\beta(s)}{s^{1+c}}ds\mbox{ for all }t\geq t_0\right\}
\end{eqnarray*}
and assume that
\begin{eqnarray*}
\alpha(t)\leq\frac{1}{t^{1-c}}\int_{t_0}^{t}f(s,u(s))ds\leq\beta(t),\qquad u\in F_{a,b,c},\thinspace t\geq t_0,
\end{eqnarray*}
where the function $f:[t_0,+\infty)\times{\R}\rightarrow[0,+\infty)$ is continuous. As a plus,
\begin{eqnarray*}
\vert f(t,u_{2}(t))-f(t,u_{1}(t))\vert\leq\frac{k(t)}{t}\vert u_{2}(t)-u_{1}(t)\vert
\end{eqnarray*}
for all $u_{1,2}\in F_{a,b,c}$ and $t\geq t_0$. Here, the function $k:[t_0,+\infty)\rightarrow[0,+\infty)$ is continuous and such that
\begin{eqnarray}
\varsigma=\frac{1}{c}\int_{t_0}^{+\infty}k(t)dt<1.\label{estimmustafa5}
\end{eqnarray}

Then, the boundary value problem
\begin{eqnarray}
\left\{
\begin{array}{ll}
x^{\prime\prime}+f(t,x)=0,\qquad t\geq t_0,\\
x(t)\geq b,\qquad t\geq t_0\\
x(t)=at+O(t^{1-c})\qquad\mbox{when }t\rightarrow+\infty
\end{array}
\right.\label{MustRogoPDE1}
\end{eqnarray}
has a unique solution $x(t)$ with the property that
\begin{eqnarray}
\alpha(t)\leq t^{c}\left[\frac{x(t)-b}{t}-x^{\prime}(t)\right]\leq\beta(t),\qquad t\geq t_0.\label{MustRogoPDE2}
\end{eqnarray}
\end{teorema}

\textbf{Proof.} Introduce the set $G$ by the formula
\begin{eqnarray*}
G=\{v\in C([t_0,+\infty),{\R}\}:-t^{-c}\beta(t)\leq v(t)\leq-t^{-c}\alpha(t)\mbox{ for all }t\geq t_0\}.
\end{eqnarray*}
The distance between the elements $v_1$ and $v_2$ of the set $D$ has the formula
\begin{eqnarray*}
d(v_1,v_2)=\sup\limits_{t\geq t_0}\{t^{c}\vert v_{1}(t)-v_2(t)\vert\}
\end{eqnarray*}
and the metric space $S=(G,d)$ is complete.

We define the operator $T:G\rightarrow C([t_0,+\infty),{\R})$ taking into account the comments in \cite[Theorem 8]{Mustafaglasgow}. Precisely,
\begin{eqnarray*}
T(v)(t)=-\frac{1}{t}\int_{t_0}^{t}sf\left(s,as+b-s\int_{s}^{+\infty}\frac{v(\tau)}{\tau}d\tau\right)ds,\qquad v\in G,\thinspace t\geq t_0.
\end{eqnarray*}
It is easy to notice that $T(G)\subseteq G$, since the mapping $t\mapsto at+b-t\int_{t}^{+\infty}\frac{v(s)}{s}ds$ belongs to $F_{a,b,c}$ for
all $v\in G$.

The estimates given by
\begin{eqnarray*}
&&t^{c}\vert T(v_2)(t)-T(v_1)(t)\vert\leq\frac{1}{t^{1-c}}\int_{t_0}^{t}sk(s)\int_{s}^{+\infty}\frac{\vert
v_2(\tau)-v_1(\tau)\vert}{\tau}d\tau ds\\
&&\leq\frac{1}{ct^{1-c}}\int_{t_0}^{t}s^{1-c}k(s)ds\cdot d(v_1,v_2)\leq\frac{1}{c}\int_{t_0}^{t}k(s)ds\cdot d(v_1,v_2)
\end{eqnarray*}
show that \textit{the operator $T$ is a contraction of coefficient $\varsigma$ in $S$}.

By denoting with $v_0$ its fixed point, where $v_0\in G$, the solution we are looking for has the formula $x(t)\equiv
at+b-t\int_{t}^{+\infty}\frac{v_0(s)}{s}ds$. $\square$

\begin{observatia}
\emph{We notice that, in the circumstances of Theorem \ref{MustRogoPDES}, if we assume that $\lim\limits_{t\rightarrow+\infty}t^{1-c}\alpha(t)$
$=\lim\limits_{t\rightarrow+\infty}t^{1-c}\beta(t)=d\in[0,+\infty)$, all the elements of the set $F_{a,b,c}$, and, in particular, the solution
of problem (\ref{MustRogoPDE1}), (\ref{MustRogoPDE2}), have the asymptotic profile $u(t)=At+B+o(1)$ as $t\rightarrow+\infty$, where $A=a$ and
$B=b+d$. Moreover, $\int_{t_0}^{+\infty}sf(s,u(s))ds=d$ for all $u\in F_{a,b,c}$.}
\end{observatia}

\begin{corolarul}\label{corMustRogoPDEs}
Set $t_0,\lambda\geq1$, $a,b\geq0$, $c\in(0,1]$ and $\varepsilon\in(0,1)$. Assume that the continuous function
$q:[t_0,+\infty)\rightarrow[0,+\infty)$ satisfies the conditions
\begin{eqnarray*}
\lambda(a+\varepsilon)^{\lambda-1}I_{c}<c\qquad\mbox{and}\qquad\frac{b}{t_0}+(a+\varepsilon)^{\lambda}\frac{I_c}{ct_{0}^{c}}<\varepsilon,
\end{eqnarray*}
where $I_c=\int_{t_0}^{+\infty}t^{\lambda+c}q(t)dt$.

Then, the equation (\ref{ecEF1}) admits the solution $x:[t_0,+\infty)\rightarrow[b,+\infty)$ with the asymptotic profile $x(t)=at+O(t^{1-c})$ as
$t\rightarrow+\infty$ for which
\begin{eqnarray*}
a^{\lambda}\cdot\frac{1}{t}\int_{t_0}^{t}s^{\lambda+1}q(s)ds\leq
\frac{x(t)-b}{t}-x^{\prime}(t)\leq(a+\varepsilon)^{\lambda}\cdot\frac{1}{t^{c}}\int_{t_0}^{t}s^{\lambda+c}q(s)ds,\quad t\geq t_0.
\end{eqnarray*}

In particular, $\frac{x(t)-b}{t}-x^{\prime}(t)=O(t^{-c})$ when $t\rightarrow+\infty$.
\end{corolarul}

\textbf{Proof.} Introduce the functions
\begin{eqnarray*}
\alpha(t)=\frac{a^{\lambda}}{t^{1-c}}\int_{t_0}^{t}s^{\lambda+1}q(s)ds,\qquad\beta(t)=(a+\varepsilon)^{\lambda}\int_{t_0}^{t}s^{\lambda+c}q(s)ds,\qquad
t\geq t_0.
\end{eqnarray*}
Obviously, $\beta(t)\leq(a+\varepsilon)^{\lambda}I_c$ in $[t_0,+\infty)$.

We have that
\begin{eqnarray*}
&&\alpha(t)\leq\frac{1}{t^{1-c}}\int_{t_0}^{t}sq(s)\left[as+b-s\int_{s}^{+\infty}\frac{v(\tau)}{\tau}d\tau\right]^{\lambda}ds\\
&&\leq\frac{1}{t^{1-c}}\int_{t_0}^{t}sq(s)\left(as+b+s\int_{s}^{+\infty}\frac{\Vert\beta\Vert_{\infty}}{{\tau}^{1+c}}d\tau\right)^{\lambda}ds\\
&&\leq\frac{1}{t^{1-c}}\int_{t_0}^{t}sq(s)\left[\left(a+\frac{b}{t_0}+\frac{(a+\varepsilon)^{\lambda}I_c}{ct_{0}^{c}}\right)s\right]^{\lambda}ds\\
&&\leq\frac{1}{t^{1-c}}\int_{t_0}^{t}s^{\lambda+1}q(s)(a+\varepsilon)^{\lambda}ds\leq\beta(t),\qquad v\in D,\thinspace t\geq t_0,
\end{eqnarray*}
and respectively
\begin{eqnarray*}
&&\vert f(t,u_2(t))-f(t,u_1(t))\vert\leq\lambda
t^{\lambda}q(t)\left\{\frac{1}{t}\left[at+b+t\int_{t}^{+\infty}\frac{\beta(s)}{s^{1+c}}ds\right]\right\}^{\lambda-1}\\
&&\times\frac{\vert u_2(t)-u_1(t)\vert}{t}\leq\lambda t^{\lambda}q(t)(a+\varepsilon)^{\lambda-1}\frac{\vert
u_2(t)-u_1(t)\vert}{t}=\frac{k(t)}{t}\vert u_2(t)-u_1(t)\vert
\end{eqnarray*}
for all $u_{1,2}\in F_{a,b,c}$ and $t\geq t_0$. $\square$

We have explained the significance of such results in the Introduction. The solution $x(t)$ of equation (\ref{ecEF1}) described in Corollary
\ref{corMustRogoPDEs} verifies, evidently, the inequality
\begin{eqnarray}
x^{\prime}(t)<\frac{x(t)}{t},\qquad t>t_0.\label{semnnegpseudowronsk}
\end{eqnarray}
On the other hand, such an estimate is not immediate for the solution $x(t)$ given at Corollary \ref{cormustafa1}. More precisely, we have
\begin{eqnarray*}
\frac{x(t)}{t}&=&A+\frac{1}{t}\int_{t_0}^{t}\int_{s}^{+\infty}q(\tau)[x(\tau)]^{\lambda}d\tau
ds=A+\int_{t}^{+\infty}q(\tau)[x(\tau)]^{\lambda}d\tau\\
&+&\left\{\frac{1}{t}\int_{t_0}^{t}sq(s)[x(s)]^{\lambda}ds
-\frac{t_0}{t}\int_{t_0}^{+\infty}q(\tau)[x(\tau)]^{\lambda}d\tau\right\}\\
&=&x^{\prime}(t)+o(1)\qquad\mbox{when }t\rightarrow+\infty,
\end{eqnarray*}
cf. \cite[Eq. (42)]{Agarwal_et_al}.

A general theorem, of Hale-Onuchic type, regarding the existence of solutions to nonlinear differential equations that are subjected to the
inequality (\ref{semnnegpseudowronsk}) can be read in \cite[Theorem 5]{Agarwal_et_al}. The next result is an adaptation of it to the
circumstances of equation (\ref{ecEF1}).

\begin{teorema}\label{MustafaRogovPDESapplic}Fix $t_0,\lambda\geq1$ and $c\geq0$, $d>0$
such that
\begin{eqnarray*}
\max\left\{\lambda(c+d)^{\lambda-1},\frac{(c+d)^{\lambda}}{d}\right\}\cdot\int_{t_0}^{+\infty}t^{\lambda}q(t)dt<1,
\end{eqnarray*}
where the function $q:[t_0,+\infty)\rightarrow[0,+\infty)$ is continuous and has (eventually) isolated zeros.

Then, the equation (\ref{ecEF1}) possesses the solution $x(t)$ defined in $[t_0,+\infty)$ with the property that
\begin{eqnarray*}
c-d\leq x^{\prime}(t)<\frac{x(t)}{t}\leq c+d,\qquad t>t_0.
\end{eqnarray*}
As a plus, the solution has the asymptotic profile $x(t)=ct+o(t)$ when $t\rightarrow+\infty$.
\end{teorema}

\textbf{Proof.} Consider $S=(D,d)$ the metric space given by the formulas
\begin{eqnarray*}
D=\{u\in C([t_0,+\infty),{\R}):ct\leq u(t)\leq (c+d)t\mbox{ for all }t\geq t_0\}
\end{eqnarray*}
and
\begin{eqnarray*}
d(u_1,u_2)=\sup_{t\geq t_0}\left\{\frac{\vert u_1(t)-u_2(t)\vert}{t}\right\},\qquad u_{1,2}\in D.
\end{eqnarray*}

For the operator $T:D\rightarrow C([t_0,+\infty),{\R})$ defined by
\begin{eqnarray*}
T(u)(t)=t\left\{c+\int_{t}^{+\infty}\frac{1}{s^{2}}\int_{t_0}^{s}\tau q(\tau)[u(\tau)]^{\lambda}d\tau ds\right\},\qquad u\in D,\thinspace t\geq
t_0,
\end{eqnarray*}
the next estimates hold
\begin{eqnarray*}
c&\leq&\frac{T(u)(t)}{t}= c+\int_{t}^{+\infty}\frac{1}{s^{2}}\int_{t_0}^{s}\tau^{\lambda+1}
q(\tau)\left[\frac{u(\tau)}{\tau}\right]^{\lambda}d\tau ds\\
&\leq& c+(c+d)^{\lambda}\int_{t}^{+\infty}\frac{1}{s^{2}}\int_{t_0}^{s}\tau^{\lambda+1}
q(\tau)d\tau ds\\
&\leq&c+(c+d)^{\lambda}\left[\frac{1}{t}\int_{t_0}^{t}\tau^{\lambda+1}q(\tau)d\tau+\int_{t}^{+\infty}\tau^{\lambda}q(\tau)d\tau\right]\\
&\leq&c+(c+d)^{\lambda}\int_{t_0}^{+\infty}\tau^{\lambda}q(\tau)d\tau<c+d
\end{eqnarray*}
and
\begin{eqnarray*}
&&\frac{\vert
T(u_2)(t)-T(u_{1})(t)\vert}{t}\\
&&\leq\int_{t}^{+\infty}\frac{1}{s^{2}}\int_{t_0}^{s}\tau^{\lambda+1}q(\tau)
\left\vert\left(\frac{u_{2}(\tau)}{\tau}\right)^{\lambda}-\left(\frac{u_{1}(\tau)}{\tau}\right)^{\lambda}\right\vert
d\tau ds\\
&&\leq\int_{t}^{+\infty}\frac{1}{s^{2}}\int_{t_0}^{s}\tau^{\lambda+1}q(\tau)\left[\lambda(c+d)^{\lambda-1}\right]
\frac{\vert u_2(\tau)-u_1(\tau)\vert}{\tau}d\tau ds\\
&&\leq\lambda(c+d)^{\lambda-1}\left[\frac{1}{t}\int_{t_0}^{t}\tau^{\lambda+1}q(\tau)d\tau+\int_{t}^{+\infty}\tau^{\lambda}q(\tau)d\tau\right]
d(u_1,u_2)\\
&&\leq\lambda(c+d)^{\lambda-1}\int_{t_0}^{+\infty}\tau^{\lambda}q(\tau)d\tau\cdot d(u_1,u_2)=\vartheta\cdot d(u_1,u_2).
\end{eqnarray*}
These imply that $T(D)\subseteq D$, respectively $T:S\rightarrow S$ is a contraction of coefficient $\vartheta$.

By denoting with $x$, where $x\in D$, the fixed point of operator $T$, we notice that
\begin{eqnarray*}
x^{\prime}(t)=[T(x)]^{\prime}(t)=\frac{x(t)}{t}-\frac{1}{t}\int_{t_0}^{t}\tau q(\tau)[x(\tau)]^{\lambda}d\tau<\frac{x(t)}{t}
\end{eqnarray*}
and
\begin{eqnarray*}
x^{\prime}(t)\geq c-\frac{1}{t}\int_{t_0}^{t}\tau q(\tau)[x(\tau)]^{\lambda}d\tau\geq c-(c+d)^{\lambda}\int_{t_0}^{t}\tau^{\lambda}q(\tau)d\tau
\end{eqnarray*}
for all $t>t_0$. $\square$

\section{An application to equation (\ref{main_eq})}

To give an application of Theorem \ref{MustafaRogovPDESapplic}, assume that, in accordance with \cite{Constantin1997,Ehrnstrom,NS1980}, the
functions $f:\overline{G}_A\times {\R}\rightarrow {\R}$ and $g:[A,+\infty)\rightarrow[0,+\infty)$ are locally H\"{o}lder continuous. Moreover,
\begin{eqnarray*}
0\leq f(x,u)\leq a(\vert x\vert)u,\quad x\in \overline{G}_A,\,u\in [0,\varepsilon],
\end{eqnarray*}
for a certain $\varepsilon >0$. Here, the function $a:[A,+\infty)\rightarrow[0,+\infty)$ is continuous and such that
\begin{eqnarray*}
\int^{+\infty}ta(t)dt<+\infty.
\end{eqnarray*}

Following the presentations in \cite{Agarwal_et_al,MustafaPDEsJMAA,MustafaRogovPDEs}, if $u_2(x)$ is a positive radially symmetric solution of
the linear elliptic equation
\begin{eqnarray}
\Delta u+a(\vert x\vert)u=0,\qquad\vert x\vert>A,\label{PDE_intermed1}
\end{eqnarray}
such that $x\cdot\nabla u_2(x)\leq0$ in $G_A$ and $u_1(x)$ is a nonnegative radially symmetric solution of the linear elliptic equation
\begin{eqnarray}
\Delta u+g(\vert x\vert)x\cdot\nabla u=0,\qquad\vert x\vert>A,\label{PDE_intermed2}
\end{eqnarray}
that satisfies the inequality $u_{1}(x)\leq u_2(x)$ throughout $G_A$, the equation (\ref{main_eq}) will possess a solution $u(x)$, not
necessarily with radial symmetry, such that
\begin{eqnarray*}
u_{1}(x)\leq u(x)\leq u_{2}(x),\qquad \vert x\vert>A.
\end{eqnarray*}

We introduce the quantities $u_{1,2}(x)=\frac{h_{1,2}(s)}{s}$, where
\begin{eqnarray*}
\vert x\vert=\left(\frac{s}{n-2}\right)^{\frac{1}{n-2}}=\beta(s).
\end{eqnarray*}

Now, the existence of solution $u_2$ to the equation (\ref{PDE_intermed1}) is implied by the (eventual) existence of a solution $h_{2}(s)$ of
the equation
\begin{eqnarray*}
h^{\prime\prime}+\frac{\beta(s)\beta^{\prime}(s)}{(n-2)s}a(\beta(s))h=0,\qquad s\geq s_0\geq1,\quad (\mbox{here, }\beta(s_0)>A)
\end{eqnarray*}
such that, in $[s_0,+\infty)$,
\begin{eqnarray}
\rho C\leq h^{\prime}(s)<\frac{h(s)}{s}\leq C\qquad\mbox{for given }C\in(0,\varepsilon),\thinspace\rho\in(0,1).\label{PDE_intermed3}
\end{eqnarray}
Since
\begin{eqnarray*}
\int_{s_0}^{+\infty}s\left[\frac{\beta(s)\beta^{\prime}(s)}{(n-2)s}a(\beta(s))\right]ds =\frac{1}{n-2}\int_{\beta(s_0)}^{+\infty}\tau
a(\tau)d\tau<+\infty,
\end{eqnarray*}
the hypotheses of Theorem \ref{MustafaRogovPDESapplic} are verified. So, there exists the supersolution $u_2(x)$ of equation (\ref{main_eq}).

Further, the problem of existence for the subsolution $u_1(x)$ which satisfies the equation (\ref{PDE_intermed2}) reduces to the existence of a
nonnegative solution $h_{1}(s)$ to the equation
\begin{eqnarray*}
h^{\prime\prime}+k(s)\left(h^{\prime}-\frac{h}{s}\right)=0,\qquad s\geq s_0,
\end{eqnarray*}
where $k(s)\equiv\beta(s)\beta^{\prime}(s)g(\beta(s))$ is a continuous nonnegative-valued function. By fixing $h_0\in(0,s_0\rho C)$ -- see
(\ref{PDE_intermed3}), we have
\begin{eqnarray*}
h_{1}(s)=s\left(\frac{h_0}{s_0}+\int_{s_0}^{s}\frac{H(\tau)}{\tau^{2}}d\tau\right),\quad H(\tau)=-\exp\left(-\int_{s_0}^{\tau}k(\xi)d\xi\right)
\end{eqnarray*}
for all $s\geq\tau\geq s_0$. In this way,
\begin{eqnarray*}
\frac{h_0-1}{s_0}\leq\frac{h_1(s)}{s}\leq\frac{h_0}{s_0},\qquad s\geq s_0.
\end{eqnarray*}

In conclusion, we have demonstrated that the equation (\ref{main_eq}) admits a bounded solution $u$ estimated by
\begin{eqnarray*}
\frac{h_0-1}{s_0}\leq u(x)\leq C,\qquad x\in G_{\beta(s_0)}.
\end{eqnarray*}

The result improves and clarifies the inferences of \cite[Section 3]{MustafaRogovPDEs}.

\textbf{Acknowledgement} The second author was financed by the Romanian AT Grant 30C/21.05.2007 with the CNCSIS code 100. Both authors are grateful to a referee for several useful comments that have improved the presentation.

\end{document}